\documentclass[12pt,draftcls,onecolumn]{IEEEtran}

\usepackage{graphicx}          

\usepackage{subfig}
\usepackage{url}
\usepackage{amssymb}
\usepackage{psfrag,graphicx}
\usepackage{verbatim}
\usepackage{amsmath,epsfig}
\usepackage{cite}
\usepackage{amssymb}
\usepackage{afterpage}
\usepackage{algorithm}
\usepackage{algorithmic}
\usepackage{threeparttable}

\newcommand{\ud}{\mathrm{d}}
\newcommand{\E}{\mathbf{E}}

\title{Iterative Source--Channel Coding Approach to Witsenhausen's Counterexample
}

\author{Johannes Kron, Ather Gattami, Tobias J. Oechtering, Mikael Skoglund\\ 
School of Electrical Engineering \\
Royal Institute of Technology KTH\\
Stockholm, Sweden}

\begin{document}
\maketitle

\begin{abstract}                          
In 1968, Witsenhausen introduced his famous counterexample where
he showed that even in the simple linear quadratic static team
decision problem, complex nonlinear decisions could outperform any
given linear decision. This problem has served as a benchmark
problem for decades where researchers try to achieve the optimal
solution. This paper introduces a systematic iterative
source--channel coding approach to solve problems of the
Witsenhausen Counterexample-character. The advantage of the
presented approach is its simplicity. Also, no assumptions are
made about the shape of the space of policies. The minimal cost
obtained using the introduced method is $\mathbf{0.16692462}$, which
is the lowest known to date.
\end{abstract}


\section{Introduction}
The most fundamental problem in control theory, name\-ly the static output feedback problem has been open since
the birth of control theory. The question is whether there is an efficient algorithm that can decide existence and find stabilizing controllers, linear or nonlinear, based on imperfect measurements and given memory.
The static output feedback problem is just an instance of the problem of control with information structures
imposed on the controllers, which has been very challenging for decision theory
researchers. In 1968, Witsenhausen \cite{witsenhausen:1968} introduced his famous counterexample:
\begin{align}
    \inf_{\gamma_1(\cdot), \gamma_2(\cdot)} \E\hspace{1mm} [k^2 \gamma_1^2(X_0) +
 X_2^2 ] \label{eq:min}
\end{align}
where
\begin{align}
    X_1 &= \gamma_1(X_0) + X_0, \\
    X_2 &=  X_1 - \gamma_2(Y_2),
\end{align}
\begin{align}
    Y_1 &= X_0, \\
    Y_2 &= X_1 + W,
\end{align}
$X_0 \sim  N(0, \sigma^2)$, and $W \sim  N(0, 1)$. Here we have
two decision makers, one corresponding to $\gamma_1$ and the other
to $\gamma_2$. The problem is a two-stage linear quadratic Gaussian control problem,
where the cost at the first time-step is $ \E[k^2 \gamma_1^2(X_0)]$ and $\E [X_2^2]$ at
the second one. At the first time-step, the controller has full state measurement, $Y_1 = X_0$.
At the second time-step, it has imperfect state measurement, $Y_2 = X_1 + W$.
What is different to the classical output feedback problem, is that the controller at the second stage does not have information from the past since it has no information about the output $Y_1$. Thus, the controller is restricted to
be a static output feedback controller. Witsenhausen showed that even in the simple linear
quadratic Gaussian control problem above, complex nonlinear
decisions could outperform any given linear decision. This problem
has served as a benchmark problem for decades where researchers
try to achieve the optimal solution. It has been pointed out that
the problem is complicated due to a so called ``signaling-incentive'', where decisions are not only chosen to
minimize a given cost, but also to encode information in the
decisions in order to signal information to other decision makers
in the team. In the example above, decision maker 2 measures $Y_2
= X_0+\gamma_1(X_0)+W$, so its measurement is affected by decision
maker 1 through $\gamma_1$. Hence, decision maker 1 not only tries
to optimize the quadratic cost in (\ref{eq:min}), but also
\textit{signal} information about $X_0$ to decision maker 2
through its decision, $\gamma_1(X_0)$.

\begin{figure}[t]
    \center
    \psfrag{X0}[][]{\small{$X_0$}}
    \psfrag{W}[][]{\small{$W$}}
    \psfrag{X1}[][]{\small{$X_1$}}
    \psfrag{Y2}[][]{\small{$Y_2$}}
    \psfrag{s1}[][]{\small{$\gamma_1$}}
    \psfrag{s2}[][]{\small{$\gamma_2$}}
    \includegraphics[width = .99\columnwidth]{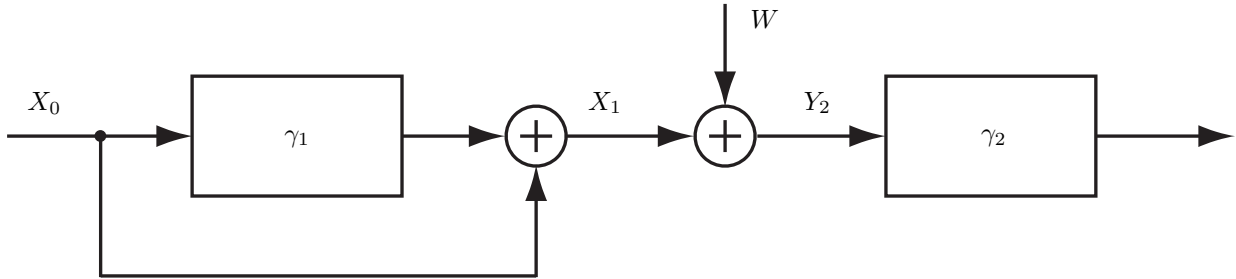}
    \caption{Schematic view of the system.}
    \label{fig:system}
\end{figure}


Previous work has been pursued on understanding the Witsenhausen
Counterexample. Suboptimal solutions where found in \cite{ho:chu}
studied variations of the problem when the signaling incentive was
eliminated. In \cite{ho:1978, ho:1980}, connections to information
theory where studied. An extensive study of the information
theoretic connection was made in \cite{bansal:basar}, where it was
shown that coupling between decision makers in the cost function
introduced the nonlinear behavior of the optimal strategies. An
ordinal optimization approach was introduced in
\cite{Deng:Ho:1999} and a hierarchical search approach was
introduced in \cite{Lee:Lau:Ho}, where both rely on a given
structure of the decisions. The first method that showed that
optimal strategies may have ``slopes'' to the quantizations was
given in \cite{baglietto01}. Solutions with bounds are studied in
\cite{Grover:Park:Sahai}. A potential games approach in the paper
by \cite{Li09} found the best known value to the date of its
publication, namely $0.1670790$.

In this paper, we will introduce a generic method of iterative
optimization based on ideas from source--channel coding
\cite{Gadkari96, Fuldseth97, Wernersson09a, Karlsson10}, that
could be used to solve problems of the Witsenhausen Counterexample
character. The numerical solution we obtain for the benchmark
problem is of high accuracy and renders the lowest value known
to date, $\mathbf{0.16692462}$.
In the following, $p(\cdot)$ and $p(\cdot | \cdot)$ denote probability density
functions (pdfs) and conditional pdfs, respectively.

\section{Iterative Optimization}
\label{sec:opti} We will now present an iterative design algorithm, based on person-by-person optimality,
for solving the minimization in equation~(\ref{eq:min}). The method we
propose is related to the Lloyd--Max algorithm~\cite{Lloyd82, Max60, Gersho92} that is successfully used when designing quantizers. A quantizer can be described by its partition cells and their corresponding reproduction value. The partition cells define to which codeword analog values are encoded and the reproduction values define how the analog value is reproduced from the codeword. In general, there is no explicit, closed-form solution to the problem of finding the optimal quantizer~\cite{Gersho92}. The key idea of the Lloyd--Max algorithm is to assume that either the partition cells or the reproduction values are fixed; with one part fixed, it is straightforward to derive an optimal expression for its counterpart. Next one part at a time is optimized in an iterative fashion. The Lloyd--Max algorithm has been generalized and used in various joint source--channel coding applications. See for example~\cite{Farvardin87, Zeger88, Farvardin91}, where quantization for noisy channels is studied, \cite{Fuldseth97}, where bandwidth compression mappings are designed, and \cite{Wernersson09a, Karlsson10}, where systems for distributed source coding and cooperative transmission are optimized.
The original Lloyd--Max algorithm converges to the global optimum under certain conditions, however, when the system model gets more complicated, as in the joint source--channel coding problems, there are no such guarantees. The algorithm can be shown to converge, but the convergence point may be only locally optimal.

The above mentioned joint source--channel coding problems are all very similar in structure to the Witsenhausen counterexample. We therefore propose to use a generalization of the Lloyd--Max algorithm to this problem; the algorithm involves four key elements:
\begin{enumerate}
  \item Formulation of necessary conditions on $\gamma_1$ and $\gamma_2$ such that they are individually optimal given that $\gamma_2$ and $\gamma_1$, respectively, are fixed.
  \item Discretization of the ``channel'' space between $\gamma_1$ and $\gamma_2$ such that $X_1$ and the input to $\gamma_2$ are restricted to belong to a finite set $\mathcal{S}_L$.
  \item Iterative optimization of $\gamma_1$ and $\gamma_2$ to make sure that they, one at a time, fulfill their corresponding necessary conditions.
  \item Use of a technique called parameter relaxation that makes the solution less sensitive to the initialization.
\end{enumerate}

\subsection{Necessary Conditions on $\gamma_1$ and $\gamma_2$}
\label{ssec:necessary_g1}
Let us first define the function
    $\tilde{\gamma}_1(x_0) \triangleq \gamma_1(x_0) + x_0 = x_1$.
Without loss of generality, we will optimize with respect to
$\tilde{\gamma}_1$.
The cost we want to minimize is given by
\begin{equation}
\label{eq:cost}
J\triangleq \E[k^2 \gamma_1^2(X_0) +
    (X_1 - \gamma_2(Y_2))^2 ].
\end{equation}
By using Bayes' rule and assuming that $\gamma_2$ is fixed, we can rewrite the optimization as
\begin{align}
\inf_{\tilde{\gamma}_1} \iint \! p(y_2 | \tilde{\gamma}_1(x_0)) \, F(x_0, \tilde{\gamma}_1(x_0), \gamma_2(y_2)) \, \ud y_2 \, p(x_0) \ud x_0 \nonumber \\
\overset{(a)}{=} \int \Big[ \inf_{x_1 \in \mathbb{R}} \int p(y_2 | x_1) \, F(x_0, x_1, \gamma_2(y_2)) \, \ud y_2 \Big] p(x_0) \ud x_0 \nonumber
\end{align}
  where
\begin{align}
  F(x_0,x_1, \gamma_2(y_2)) = \Big( k^2 (x_1 - x_0)^2 + ( x_1 - \gamma_2(y_2))^2 \Big). \nonumber
\end{align}
In (a) we make use of Theorem 14.60 in \cite{Rockafellar2009}, which states that interchange of minimization and integration is possible under certain conditions\footnote{In our case the conditions are fulfilled because the integrand is continuous in $x_0$ and $x_1 = \tilde{\gamma}_1 (x_0)$ and the infimum is over the space of all measurable functions.}.
Furthermore, since the optimal value of $J$ is not $-\infty$, the theorem states that $\tilde{\gamma}_1(\cdot)$ can be defined in a pointwise manner. Consequently, a necessary condition for $\tilde{\gamma}_1$ to be optimal
is given by
\begin{align}
    \tilde{\gamma}_1(x_0) = &\arg \min_{x_1 \in \mathbb{R}} \Big( \int p(y_2 | x_1) \,\, F(x_0,x_1,\gamma_2(y_2)) \,\,  \ud y_2 \Big)\label{eq:opt_g1}
\end{align}
for almost every $x_0\in \mathbb{R}$.


If we next assume that $\gamma_1$ is
fixed, we see that the first term in (\ref{eq:cost}) is a
constant. The minimization of $J$ with respect to $\gamma_2$ is
therefore equivalent to
\begin{align}
    \inf_{\gamma_2(\cdot)} \E[(X_1 - \gamma_2(Y_2))^2],
\end{align}
which is the mean-squared error (MSE). It is well known that the MSE is minimized by the conditional expected value; hence,
\begin{align}
    \gamma_2(y_2) = \E[X_1 | y_2] \label{eq:opt_g2}
\end{align}
for almost every $y_2 \in \mathbb{R}$, is a necessary condition for $\gamma_2(y_2)$ to be optimal.

\subsection{Discretization}
\label{ssec:discrete}
The expressions given in~(\ref{eq:opt_g1}) and~(\ref{eq:opt_g2}) are impractical to use in our design algorithm because they require the functions to be specified for infinitely many input values. To get around this problem we introduce a discrete set
\begin{align}
\mathcal{S}_L = \Big\{-\Delta \frac{L-1}{2}, -\Delta \frac{L-3}{2},
  \ldots, \Delta \frac{L-3}{2}, \Delta \frac{L-1}{2} \Big\}, \label{eq:S}
\end{align}
where $L \in \mathbb{N}$ and $\Delta \in \mathbb{R}_+$ are two
parameters that determine the number of points and the spacing
between the points, respectively. Next, we impose the constraint
$x_1 \in \mathcal{S}_L$, that is, the output of $\tilde{\gamma}_1$
can only take one out of a finite number of values. In a similar
way, the input to $\gamma_2$ is discretized such that,
\begin{align}
    \gamma_2(y_2) = \tilde{\gamma}_2 ( \tilde{y}_2 ), \quad \tilde{y}_2 = Q_{\mathcal{S}_L}(y_2)\in \mathcal{S}_L,
\end{align}
where $Q_{\mathcal{S}_L}(y_2)$ maps $y_2$ to the closest point in the set $\mathcal{S}_L$. $\gamma_2$ can now be stored in the form of a lookup table where each point in $\mathcal{S}_L$ is associated with an output value. The approximation of the real space with $\mathcal{S}_L$ can be made more and more accurate by decreasing $\Delta$ and increasing $L$\footnote{While decreasing $\Delta$, one has to increase $L$ to make sure that $\max (x \in \mathcal{S}_L) = \Delta (L-1)/2$ does not decrease.}. Finally, since $X_0$ is still infinite-dimensional, we use Monte-Carlo samples of $X_0$ to represent the input to $\tilde{\gamma}_1$.
$\tilde{\gamma}_1$ is now specified by evaluating
\begin{align}
    \tilde{\gamma}_1 (x_0) = \arg \min_{x_1 \in \mathcal{S}_L} \sum_{\tilde{y}_2 \in \mathcal{S}_L} p(\tilde{y}_2 | x_1) \,\, F(x_0,x_1,\tilde{\gamma}_2(\tilde{y}_2))\,\, \label{eq:opt_g1d}
\end{align}
for each of the Monte-Carlo samples that represent $X_0$. In a similar way, $\tilde{\gamma}_2$ can be expressed as
\begin{align}
    \tilde{\gamma}_2(\tilde{y}_2) = \E[X_1 | \tilde{y}_2] \label{eq:opt_g2d},
\end{align}
for all $\tilde{y}_2 \in \mathcal{S}_L$, where the expectation with respect to $X_0$ is evaluated by using the Monte-Carlo samples.

\subsection{Design Algorithm Using Parameter Relaxation}
\label{ssec:design}
Given the above expressions for $\tilde{\gamma}_1$ and $\tilde{\gamma}_2$
it is possible to optimize the system iteratively. One
common problem with iterative techniques is that the final solution will depend on the initialization
of the algorithm. If the initialization is bad we are likely to
end up in a poor local minimum. One method that has proven to be helpful in counteracting this in joint source--channel coding is noisy channel relaxation (NCR) \cite{Gadkari96, Fuldseth97,
Wernersson09a, Karlsson10}. In this paper, we use a generalization of NCR which we call parameter relaxation (PR). The idea of PR is to first define a parameter space $\mathcal{P}$ that include relevant system parameters such as noise variance, power constraints, etc. Assuming that we have found a system that performs well for a system parameter $\eta_n \in \mathcal{P}$, this system is then used as initialization when designing a new system for a parameter $\eta_{n+1} = \eta_n + \epsilon \in \mathcal{P}$. This update procedure is continued until $\eta_{n} = \eta_T$, which is the target system parameter (i.e., the system parameter for which we want to find the optimized system.).

The problem in the PR method is to determine a good starting point $\eta_0$ as well as the path to reach $\eta_T$. In joint source--channel coding, the most common parameter to change is the noise variance of the channel. The optimization starts with a high noise variance which is gradually decreased to the target noise variance, hence the name NCR. In the Witsenhausen setup, we have found that the parameter $k$ is useful to include in the parameter space: Design a system for a high value of $k$ first and
then gradually decrease $k$ until the desired value of $k:T$ is
reached. The reason to start with a high value of $k$ is
that the design algorithm will find a solution where $\tilde{\gamma}_1(x_0) \approx x_0$ in this case (i.e.,
$\gamma_1(x_0) \approx 0$) independently of $\gamma_2$. The design procedure including the PR part
is given in Algorithm 1.
Each update on line $7$ and $8$ in Algorithm 1 will decrease the cost. Since the cost is lower bounded, it is clear that the algorithm will converge. It may happen that the algorithm converges to a local optimum, however, as will be seen in the following section the \emph{local} optima we obtain are still better than any previously reported results.

\algsetup{indent=1em}
\begin{algorithm}[h!]
\caption{Design Algorithm\label{alg:design}}
\begin{algorithmic}[1]
\medskip
\REQUIRE Initial mapping of $\tilde{\gamma}_2$, the value $k_T$ for which the system should be optimized and the threshold $\delta$ that determines when to stop the iterations.
\ENSURE Locally optimized $\tilde{\gamma}_1$ and $\tilde{\gamma}_2$ .
\medskip

\STATE Let $k > k_T$.
\WHILE{$k > k_T$}
\STATE Decrease $k$ according to some scheme (e.g., linearly).
\STATE Set the iteration index $i=0$ and $J^{(0)}=\infty$.
\REPEAT
\STATE Set $i = i + 1$
\STATE Find the optimal $\tilde{\gamma}_1$ by using \eqref{eq:opt_g1d}.
\STATE Find the optimal $\tilde{\gamma}_2$ by using \eqref{eq:opt_g2d}.
\STATE Evaluate the cost function $J^{(i)}$ according to \eqref{eq:cost}.

\UNTIL{$(J^{(i-1)}- J^{(i)}) / J^{(i-1)} < \delta$}
\ENDWHILE
\end{algorithmic}
\end{algorithm}

\section{Results}
\label{sec:results}
\subsection{Implementation Aspects}
For the evaluation of the design algorithm we have initially used
$L=201$ levels and chosen $\Delta(L) = 10 \sigma / (L-1)$. We have
used $400000$ Monte-Carlo samples in the final optimizations to
represent $X_0$. Since it is known that the optimal $\gamma_1$ is
symmetric about origin \cite{witsenhausen:1968}, we have restricted
$\tilde{\gamma}_1$ to have this symmetry by generating only
positive Monte-Carlo samples and thereafter reflecting the
resulting $\tilde{\gamma}_1$-function for negative values of
$x_0$.
To be able to compare our results to previously reported results,
we have set $\sigma=5$ and $k_T=0.2$. However, since we are using the PR method, we
have initially used the value $k = 3$ and decreased it according to
the series $\{3, 2, 1.5, 1, 0.6, 0.4, 0.3, 0.2\}$. Before running the design algorithm, we require $\tilde{\gamma}_2$ to be initialized, but due to the PR this has little impact on the final solution and we have used the initialization $\tilde{\gamma}_2 \equiv 0$.

Once we have obtained the solution for $k_T=0.2$, we have increased the
precision by expanding the number of points in the discrete set
from $L$ to $L'$ and updated $\tilde{\gamma}_2$ according to
\begin{align}
    \tilde{\gamma}_2^{(L')}( \tilde{y}_2 ) = \tilde{\gamma}_2^{(L)} ( Q_{\mathcal{S}_L}(\tilde{y}_2)  )
\end{align}
for all $\tilde{y}_2 \in \mathcal{S}_{L'}$. Thereafter the inner part of the design algorithm, that is, lines 4--10, have been run again to obtain a system optimized for the increased number of points $L'$. By repeating this refinement, the precision increases and the cost decreases as will be shown later. This method of refining the precision is similar to the one-way multigrid algorithm that is analyzed in~\cite{Chow91}. The evaluations of~(\ref{eq:opt_g1d}) and~(\ref{eq:opt_g2d}) have been done using an exhaustive search, therefore, the run time is exponential in the number of levels $L$. 

\subsection{Numerical Results}
\begin{figure}[t]
    \center
    \psfrag{X}[][]{\small{$x_0$}}
    \psfrag{Y}[][]{\small{$\tilde{\gamma}_1$}}
    \includegraphics[width = .99\columnwidth, trim=0 8mm 0 8mm, clip=true]{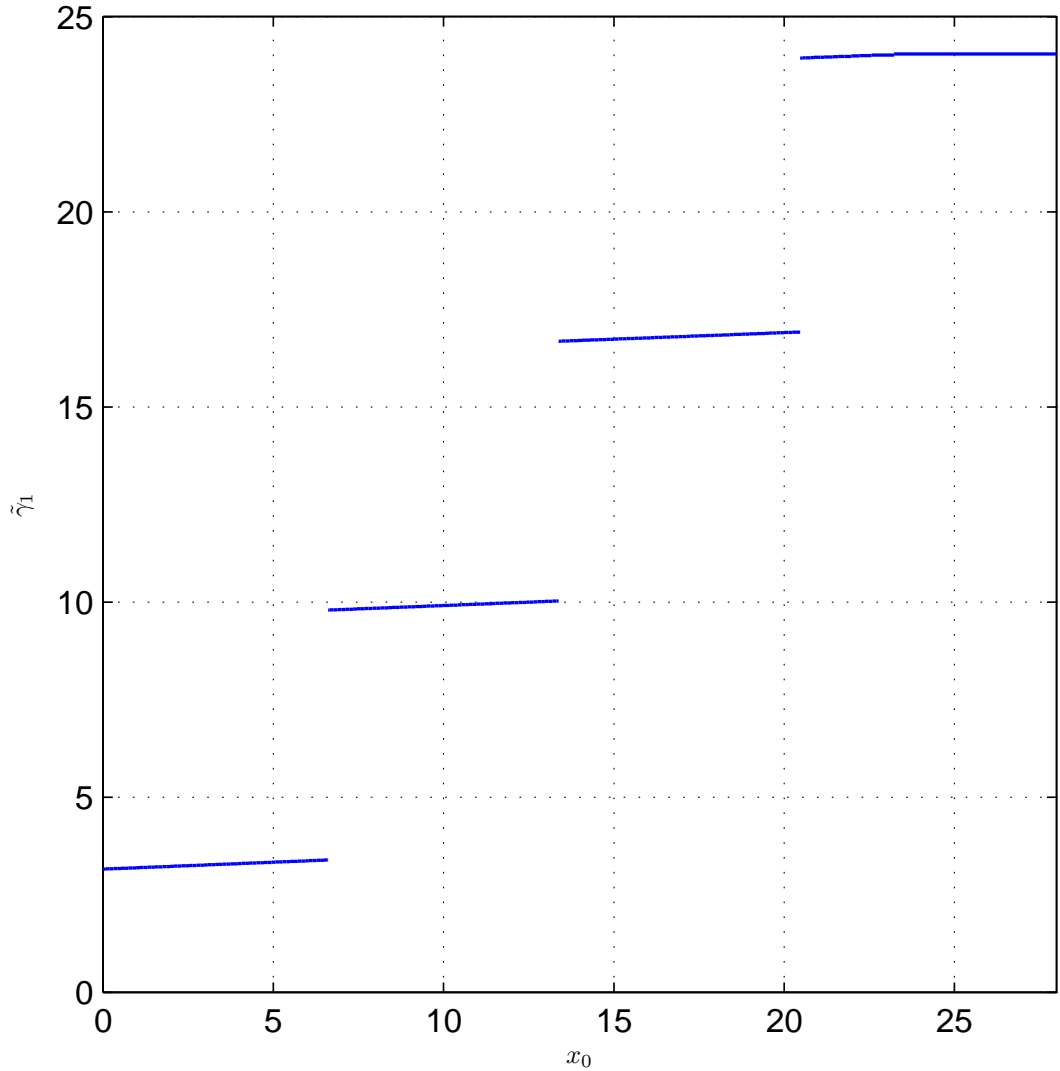}
    \caption{$4$-step solution ($L=12801$)}
    \label{fig:gamma1_4}
\end{figure}

During the first steps of the PR $k$ is high. This means that
the output of $\tilde{\gamma}_1$ should follow the input closely
to avoid large costs in the first stage. If continuous outputs
were allowed, the output would be identical to the input. However,
since we are working with a discretized system, only outputs from
the set $\mathcal{S}_L$ are feasible. As $k$ reaches $0.4$--$0.6$
the step behavior of the output appears. This particular value of $k$ where the system changes from being affine to have a more general shape is consistent with a result in~\cite{Wu11}, which states that the optimal cost is less than the optimal cost for an affine system if $k < 0.564$.
Depending on the realization of the Monte-Carlo samples we get either a $3.5$-step mapping or a $4$-step
mapping as shown in Fig.~\ref{fig:gamma1_4} (occasionally, a
$3$-step solution has occurred). The total costs for these
solutions are stated in Table~\ref{tbl:cost2}. For ease of comparison, we have also included the costs of previously reported results. As can be seen,
all our mappings have similar performance and all of them give lower costs
than the previously reported lowest cost --- $0.1670790$
\cite{Li09}.

\begin{table*}[t]%
\centering
\refstepcounter{table}
Table \arabic{table}\\
\vspace{2mm}
\subfloat[Final cost for different solutions.]{
\begin{tabular}{|r|r|c|c|c|}
  \hline
  Steps        &  Stage 1 & Stage 2 & Total Cost \\
  \hline
  Witsenhausen \cite{witsenhausen:1968}$^\dag$&  $0.40423088$ &  $0.00002232$ &  $0.40425320$  \\
  Bansal \& Bansar \cite{bansal:basar}$^\dag$ &
                  $0.36338023$ &  $0.00163460$ &  $0.36501483$  \\
  Deng \& Ho \cite{Deng:Ho:1999}$^\dag$ &  $0.13948840$ &  $0.05307290$ &  $0.19256130$  \\
  Baglietto et al. \cite{baglietto01} &         &             &  $0.1701$ \\
  Lee et al. \cite{Lee:Lau:Ho}  &  $0.13188408$ &  $0.03542912$ &  $0.16731321$  \\
  Li et al. \cite{Li09}    &    &    &  $0.1670790$ \\
  \hline
  This paper, $3$-step$^\ddag$     &  $0.13493778$ &  $0.03201113$ &  $0.16694891$  \\
  This paper, $3.5$-step$^\ddag$   &  $0.13462186$ &  $0.03230369$ &  $0.16692555$  \\
  This paper, $4$-step$^\ddag$     &  $0.13484828$ &  $0.03207634$ &  $\mathbf{0.16692462}$ \\
  \hline
\end{tabular}
\label{tbl:cost2}
}
\hspace{2mm}
\subfloat[Costs for different precisions for the $4$-step solution.]{
\begin{tabular}{|r|r|c|c|c|}
  \hline
  $L$      &    $M$    & Stage 1 & Stage 2 & Total Cost \\
  \hline
  $201$     &   $16$  & $0.121042$  & $0.057641$ &  $0.17868301$  \\
  $401$     &   $22$  & $0.130150$  & $0.038834$ &  $0.16898421$  \\
  $801$     &   $30$  & $0.135308$  & $0.032009$ &  $0.16731642$  \\
  $1601$    &   $56$  & $0.134966$  & $0.032062$ &  $0.16702853$  \\
  $3201$    &  $110$  & $0.134868$  & $0.032081$ &  $0.16694954$  \\
  $6401$    &  $210$  & $0.134859$  & $0.032071$ &  $0.16692966$  \\
  $12801$   &  $396$  & $0.134848$  & $0.032076$ &  $0.16692462$  \\
  \hline
\end{tabular}
\label{tbl:cost}
}
\vspace{-2mm}
$^\dag$ Costs obtained from \cite{Lee:Lau:Ho}. $^\ddag$ $L=12801$.
\end{table*}

In Table~\ref{tbl:cost} we show how the cost decreases as the number of points $L$ is increased. The method we use to calculate the total cost as well as some notes on the accuracy can be found in Appendix~\ref{app:totalcost}. The lowest cost we have achieved with our algorithm is $0.16692462$. The mapping that achieves this cost is the $4$-step mapping shown in Fig.~\ref{fig:gamma1_4} with $L=12801$ points.
Although the mapping contains four clear output levels it should be emphasized that each level is slightly sloped; this can be seen in Fig.~\ref{fig:gamma1_4_zoom}, where the first step has been zoomed in. It is reasonable to assume that as the precision (i.e., $L$) increases further, each step of the mapping will converge to a straight line that is slightly sloped.

\section{Comparison to Previous Results}
In this section, we will compare the presented method with 
previous methods and note some differences:

\begin{itemize}
\item No structure is assumed for the decision functions. In
\cite{Deng:Ho:1999} and  \cite{Lee:Lau:Ho}, monotonicity of the
decisions was assumed. The space of decisions is assumed to be a
normed linear space in \cite{baglietto01}.

\item The design is fully automated and little modeling needs to be done a-priori. In contrast, a significant analytic/modeling work was performed before
posing the optimization problem to be solved in
\cite{Deng:Ho:1999}, \cite{Lee:Lau:Ho}, and \cite{baglietto01}.
The first two require manual adjustments for the proper choice of
interval values and signal levels, and the third requires some
prior analysis to determine a constant ``$c$''. In \cite{Li09},
modeling work is needed in converting the problem into a potential
game.
\end{itemize}

\section{Conclusions}
In this paper, we introduced a generic method of iterative
optimization based on ideas from source--channel coding, that could
be used to solve problems of the Witsenhausen counterexample
character. The numerical solution we obtain for the benchmark
problem is of high accuracy and renders the lowest value known
to date, $\mathbf{0.16692462}$. Also, the design algorithm does not make any assumption on the
structure of the policies
--- the solutions are allowed to have arbitrary shapes (within the
restrictions imposed by the discretization). The results can
therefore be seen as a confirmation that the step-shaped behavior
is beneficial.


\appendix
\section{Calculation of the Total Cost}
\begin{figure}[t]
    \center
    \psfrag{X}[][]{\small{$x_0$}}
    \psfrag{Y}[][]{\small{$\tilde{\gamma}_1$}}
    \includegraphics[width = .7\columnwidth, trim=0 5mm 0 5mm, clip=true]{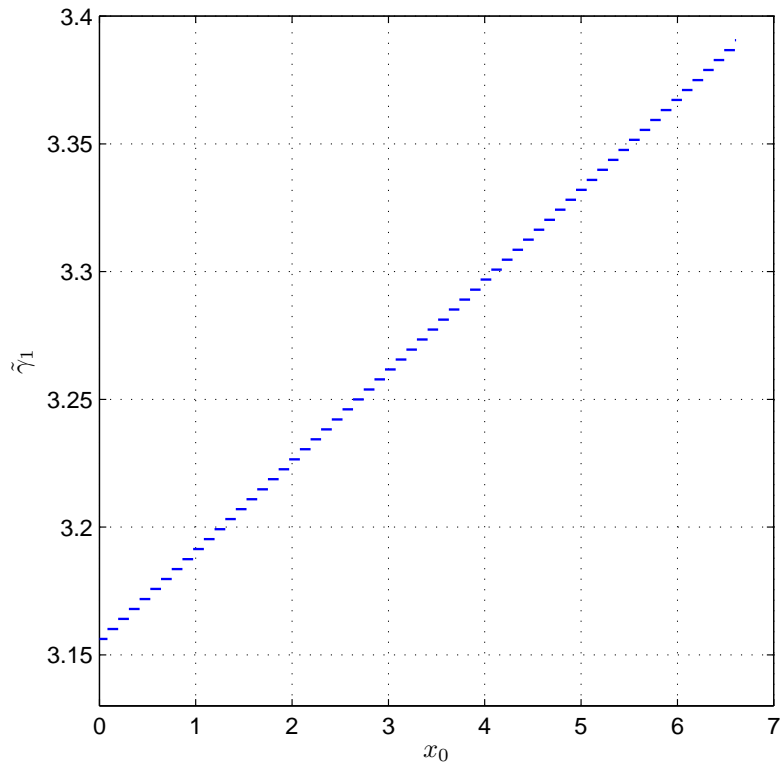}
    \caption{Detailed view of the first step in the $4$-step solution.}
    \label{fig:gamma1_4_zoom}
\end{figure}
\label{app:totalcost} In the design algorithm, $\tilde{\gamma}_1$
is specified implicitly by storing
the output symbol to which each Monte-Carlo sample is mapped. This representation is
used when evaluating the cost during the iterations in the design
algorithm. However, to evaluate the final total cost we need
higher numerical accuracy. Therefore, the first step in calculating the total
cost is to use the sample-based representation to find thresholds, $A_i$, such that $\tilde{\gamma}_1$ can be given  on the form
\begin{align}
\tilde{\gamma}_1(x_0) = \alpha_i \in \mathcal{S}_L \quad \mathrm{if}\,\, A_i \leq x_0 < A_{i+1},
\end{align}
for $i \in \{0, \ldots, M-1 \}$, with $A_0 = -\infty$ and $A_M = \infty$. That is, the sample-based representation of $\tilde{\gamma}_1$, which is explicitly defined only for the Monte-Carlo samples, is transformed to a function which is defined for all real numbers. This representation makes it possible to numerically evaluate the integrals that are needed to find the total cost
\begin{align}
   J &= \E[k^2 \gamma_1^2(X_0) + (X_1 - \gamma_2(Y_2))^2 ] \nonumber \\
     & = \underbrace{\E[k^2 (\tilde{\gamma}_1(X_0) - X_0)^2 ]}_{= J_1} + \underbrace{\E[(\tilde{\gamma}_1(X_0) - \tilde{\gamma}_2(\tilde{Y}_2))^2 ]}_{= J_2},
\end{align}
where
\begin{align}
J_1 &= \int_{x_0} p(x_0) k^2 (\tilde{\gamma}_1(x_0) - x_0)^2 \ud x_0 \nonumber \\
                &= k^2 \sum_{i=0}^{M-1} \int_{A_i}^{A_{i+1}} p(x_0) (\alpha_i - x_0)^2 \ud x_0, \\
J_2 &= \int_{x_0} \sum_{\tilde{y}_2 \in \mathcal{S}_L} p(x_0, \tilde{y}_2) (\tilde{\gamma}_1(x_0) - \tilde{\gamma}_2(\tilde{y}_2))^2  \ud x_0 \nonumber \\
&= \sum_{i=0}^{M-1} \Big\{ \!\! \sum_{\tilde{y}_2 \in
\mathcal{S}_L} \!\! P(\tilde{y}_2 | \alpha_i) (\alpha_i -
\tilde{\gamma}_2(\tilde{y}_2))^2 \Big\} \int_{A_i}^{A_{i+1}}
\!\!\!\!\!\!\! p(x_0) \ud x_0,
\end{align}
and
\begin{align}
P(\tilde{y}_2 | \alpha_i) & =
\left\{ \!\!\! \begin{array}{ll}
\displaystyle\int_{-\infty}^{\tilde{y}_2 + \Delta/2} \!\!\!\!\!\!\!\!\!\!\! p( w = y_2 - \alpha_i ) \ud y_2 & \textrm{if $\tilde{y}_2 = -\Delta \frac{L-1}{2}$}\\
\displaystyle\int_{\tilde{y}_2 - \Delta/2}^{\infty}  \!\!\!\!\!\!\!\!\!\!\! p( w = y_2 - \alpha_i ) \ud y_2 & \textrm{if $\tilde{y}_2 = \Delta \frac{L-1}{2}$}\\
\displaystyle\int_{\tilde{y}_2 - \Delta/2}^{\tilde{y}_2 + \Delta/2} \!\!\!\!\!\!\!\!\!\!\!  p( w = y_2 - \alpha_i ) \ud y_2 & \textrm{otherwise}
\end{array} \right.
\end{align}
All integrals have been calculated numerically using the Matlab
function \emph{quadl} with the tolerance specified to be $t =
10^{-18}$, which means that the absolute error of the result from quadl is not greater than $t$.
All integrands are continuous and have a smooth behavior that should cause no problem for quadl.
To upper bound the total cost, we have upper bounded each integral by adding $t$ to each individual result from quadl and reevaluated the total cost. In this way we have estimated
the absolute error to be in the order of (or less than)
$10^{-11}$. Matlab code for our calculations of the total cost,
including our decision functions can be found in
\cite{witsenhausen_code}.

%

\bibliographystyle{plain}        
\bibliography{mybib}             

\end{document}